# Fixed-point properties of the Mordukhovich differential operator

Jinlu Li

Department of Mathematics
Shawnee State University
Portsmouth, Ohio 45662 USA
jli@shawnee.edu

**Abstract**

In this paper, we investigate some fixed-point properties of the Mordukhovich differential operator of set valued mappings (or, single valued mappings) on Banach spaces. In particular, we study the fixed-point properties of the Mordukhovich differential operator for the metric projection operator onto some closed and convex subsets in Banach spaces, such as, closed balls in Banach spaces, positive cones in real spaces $l_2$ and $l_1$ and sets of polynomials in $C[0, 1]$.

In honor of Professor Dr. Sehie Park for his great contributions to Mathematics,
in particular, to Fixed Point Theory

## 1. Introduction

Let $(X, \|\cdot\|)$ be a real Banach space with topological dual space $(X^*, \|\cdot\|_*)$. Let $Y$ be a Banach spaces with topological dual space $Y^*$. If it is not confused, let $\langle\cdot, \cdot\rangle$ denote both the real canonical pairing between $X^*$ and $X$ and the real canonical pairing between $Y^*$ and $Y$. Let $\Delta$ be a nonempty subset of $X$ and let $f: \Delta \to Y$ be a single valued mapping. Traditionally, the smoothness of $f$ is described by some types of differentiability of $f$, such as Gâteaux directional differentiability, Fréchet differentiability and strict Fréchet differentiability. For example, in the case that the underlying spaces are Hilbert spaces, some differentiability and its applications of $f$ has been introduced and studied in [2, 6, 9, 10, 19]; and if the underlying spaces are Banach spaces and normed linear spaces, it has been studied in [3, 7, 12, 23].

It is well-known that set valued mappings have shown increasingly important roles in variational analysis, optimization theory, control theory, game theory, equilibrium theory, economic theory, and so forth. In contrast with the differentiability of single valued mappings, a fresh theory of generalized differentiation of set valued mappings in Banach spaces has been successfully developed. The fundamental concepts of the new generalized differentiation of set valued mappings are Mordukhovich derivatives, which are also called Mordukhovich coderivatives, or coderivatives. We briefly review the concepts of Mordukhovich coderivatives for set valued mappings in Banach spaces. For more details, for example, see [16−18].

Let $F: \Delta \rightrightarrows Y$ be a set valued mapping. For any $x \in \Delta$ and $y \in F(x)$, the Mordukhovich derivative (which is also called Mordukhovich coderivative, or coderivative) of $F$ at point $(x, y)$ is defined by (see Definitions 1.13 and 1.32 in Chapter 1 in [16])

--------------------

$$\widehat{D}^*F(x,y)(y^*) = \left\{ x^* \in X^*: \limsup_{\substack{(u,v)\to(x,y)\\ u\in\Delta \text{ and } v\in F(u)}} \frac{\langle x^*,u-x\rangle - \langle y^*,v-y\rangle}{\|u-x\|+\|v-y\|} \le 0 \right\}, \text{ for any } y^* \in Y^*. \quad (1.1)$$

The Mordukhovich derivatives have been widely applied to several branches of mathematics such as operator theory, optimization theory, approximation theory, control theory, equilibrium theory, and so forth (see [16−18]).

For any $x \in \Delta$ and $y \in F(x)$, by the above definition, $\widehat{D}^*F(x,y): Y^* \rightrightarrows X^*$ is a set valued mapping, which is called the Mordukhovich differential operator (or the Mordukhovich codifferential operator) of $F$ at the point $(x,y)$.

In particular, in the Definition (1.1), if $X = Y$, then, the Mordukhovich differential operator of $F$ at the point $(x,y)$ is $\widehat{D}^*F(x,y): X^* \rightrightarrows X^*$, which satisfies (1.1), for $y^* \in X^*$. More precisely speaking, in the case that $X = Y$, the Mordukhovich differential operator of $F$ at the point $(x,y)$ is a set valued mapping on $X^*$. Therefore, in this paper, we investigate the fixed-point properties of the mapping $\widehat{D}^*F(x,y): X^* \rightrightarrows X^*$.

In section 3, we prove some fixed-point properties of $\widehat{D}^*F(x,y)$ for general set valued mapping $F$. In section 4, we consider the special and very useful metric projection operator. Let $C$ be a nonempty closed and convex subset in $X$, and let $P_C$ denote the standard metric projection operator from $X$ onto $C$. We will use the results about the Mordukhovich derivatives of $P_C$ obtained in [13−15] to study the fixed-point properties of the Mordukhovich differential operator of the metric projection operator $P_C$, in which $C$ is a closed ball in Banach spaces, or the positive cone in real spaces $l_2$ and $l_1$, or a set of polynomials in $C[0, 1]$.

Let $C[0, 1]$ be the Banach space of all continuous real valued functions on [0, 1] with the maximum norm. Let $n$ be a positive integer and let $\mathcal{P}_n$ be the set of all polynomials with degree less than or equal to $n$, which is a subset of $C[0, 1]$. In [15], it is proved that the metric projection $P_{\mathcal{P}_n}: C[0, 1] \to \mathcal{P}_n$ is a single valued mapping. In the last subsection of Section 4, we will prove that $P_{\mathcal{P}_n}$ is a continuous mapping. By the continuity of $P_{\mathcal{P}_n}$, we provide a simplified calculation formula for fixed points of the Mordukhovich differential operator $\widehat{D}^*P_{\mathcal{P}_n}(\cdot)$ of the metric projection operator $P_{\mathcal{P}_n}$.

## 2. Preliminaries

Let $(X, \|\cdot\|)$ be a real Banach space with topological dual space $(X^*, \|\cdot\|_*)$. Let $\langle\cdot, \cdot\rangle$ denote the real canonical pairing between $X^*$ and $X$, and $\theta$, $\theta^*$ denote the origins in $X$ and $X^*$, respectively. The identity mappings on $X$ and $X^*$ are respectively denoted by $I_X$ and $I_{X^*}$. Let $\Delta$ be a nonempty subset of $X$ and let $F: \Delta \rightrightarrows X$ be a set valued mapping. The graph of $F$ is defined by the following subset in $\Delta \times X$

$$\text{gph}F = \{(x,y) \in \Delta \times X: y \in F(x)\}.$$

For $(x,y) \in \text{gph}F$, that is, for any $x \in \Delta$ and $y \in F(x)$, as a special case of (1.1), the Mordukhovich derivative of $F$ at point $(x, y)$ is defined by (This is a special case in Definitions 1.13 and 1.32 in

Chapter 1 in [16], in which $X = Y$)

$$\widehat{D}^*F(x,y)(y^*) = \left\{x^* \in X^*\colon \limsup_{\substack{(u,v)\to(x,y)\\ u\in\Delta \text{ and } v\in F(u)}} \frac{\langle x^*, u-x\rangle - \langle y^*, v-y\rangle}{\|u-x\|+\|v-y\|} \le 0\right\}$$

$$= \left\{x^* \in X^*\colon \limsup_{\substack{(u,v)\to(x,y)\\ (u,v)\in \text{gph}F}} \frac{\langle x^*, u-x\rangle - \langle y^*, v-y\rangle}{\|u-x\|+\|v-y\|} \le 0\right\}, \text{ for any } y^* \in X^*. \qquad (2.1)$$

If $(x,y) \notin \text{gph}F$, then, we define

$$\widehat{D}^*F(x,y)(y^*) = \emptyset, \text{ for any } y^* \in X^*.$$

In particular, let $f\colon \Delta \to X$ be a single valued continuous mapping. By (2.1), the Mordukhovich derivative of $f$ at point $(x, f(x))$ is defined by

$$\widehat{D}^*f(x,f(x))(y^*) \coloneqq \widehat{D}^*f(x)(y^*)$$

$$= \left\{x^* \in X^*\colon \limsup_{\substack{u\to x\\ u\in\Delta}} \frac{\langle x^*, u-x\rangle - \langle y^*, f(u)-f(x)\rangle}{\|u-x\|+\|f(u)-f(x)\|} \le 0\right\}, \text{ for any } y^* \in X^*. \qquad (2.2)$$

Let $C$ be a nonempty closed and convex subset of a Banach space $X$, and let $P_C$ denote the standard metric projection operator from $X$ onto $C$. In general, $P_C\colon X \rightrightarrows C$ is a set valued mapping. For any $x \in X$ and $y \in P_C(x)$, by (2.1), the Mordukhovich derivative of $P_C$ at point $(x, y)$ is defined by

$$\widehat{D}^*P_C(x,y)(y^*) = \left\{x^* \in X^*\colon \limsup_{\substack{(u,v)\to(x,y)\\ u\in X \text{ and } v\in P_C(u)}} \frac{\langle x^*, u-x\rangle - \langle y^*, v-y\rangle}{\|u-x\|+\|v-y\|} \le 0\right\}$$

$$= \left\{x^* \in X^*\colon \limsup_{\substack{(u,v)\to(x,y)\\ (u,v)\in \text{gph}P_C}} \frac{\langle x^*, u-x\rangle - \langle y^*, v-y\rangle}{\|u-x\|+\|v-y\|} \le 0\right\}, \text{ for any } y^* \in X^*. \qquad (2.3)$$

Some properties of Mordukhovich derivatives of $P_C$ were proved in [15] in general Banach spaces, which include $l_1$, $c$ and $C[0, 1]$.

In particular, if $X$ is a uniformly convex and uniformly smooth Banach space, then, it is well-known that $P_C\colon X \to C$ is a single valued continuous mapping. From (2.2), the Mordukhovich derivative of $P_C$ at $(x, P_C(x))$ is defined by

$$\widehat{D}^*P_C(x,P_C(x))(y^*) \coloneqq \widehat{D}^*P_C(x)(y^*)$$

$$= \left\{ x^* \in X^* : \limsup_{u \to x} \frac{\langle x^*, u-x \rangle - \langle y^*, P_C(u) - P_C(x) \rangle}{\|u-x\| + \|P_C(u) - P_C(x)\|} \le 0 \right\}, \text{ for any } y^* \in X^*. \qquad (2.4)$$

In [13−15], properties and representations of Mordukhovich derivatives of $P_C$ were provided in uniformly convex and uniformly smooth Banach spaces and Hilbert spaces, with respect to the given nonempty closed and convex subset $C$ to be closed balls, closed and convex cones and closed and convex cylinders.

## 3. Fixed point properties of Mordukhovich differentiation operator

Let $\Delta$ be a nonempty subset of a Banach space $X$. Let $F: \Delta \rightrightarrows X$ be a set valued mapping. For any $x \in \Delta$ and $y \in F(x)$, by definition (2.1), the Mordukhovich derivative of $F$ at $(x, y) \in \text{gph}F$ is $\widehat{D}^*F(x, y): X^* \rightrightarrows X^*$, which is actually a set valued operator on $X^*$. $\widehat{D}^*F(x, y)$ is called the Mordukhovich differential operator (or Mordukhovich codifferential operator) of $F$ at the point $(x, y) \in \text{gph}F$. Let $y^* \in X^*$. If

$$y^* \in \widehat{D}^*F(x, y)(y^*),$$

then, $y^*$ is a fixed point of the Mordukhovich differential operator $\widehat{D}^*F(x, y)$. The collection of all fixed points of the Mordukhovich differential operator $\widehat{D}^*F(x, y)$ is denoted by

$$\mathcal{F}\left(\widehat{D}^*F(x, y)\right) = \{y^* \in X^* : y^* \in \widehat{D}^*F(x, y)(y^*)\}.$$

By (2.1), $\mathcal{F}\left(\widehat{D}^*F(x, y)\right)$ has the following representations.

**Lemma 3.1**. *Let $\Delta$ be a nonempty subset of a Banach space $X$. Let $F: \Delta \rightrightarrows X$ be a set valued mapping. Let $x \in \Delta$ and $y \in F(x)$. Then, for any $y^* \in X^*$, we have*

$$\mathcal{F}\left(\widehat{D}^*F(x, y)\right)$$

$$= \left\{ y^* \in X^* : \limsup_{\substack{(u,v) \to (x,y) \\ (u,v) \in \text{gph}F}} \frac{\langle y^*, u-x \rangle - \langle y^*, v-y \rangle}{\|u-x\| + \|v-y\|} \le 0 \right\}$$

$$= \left\{ y^* \in X^* : \limsup_{\substack{(u,v) \to (x,y) \\ (u,v) \in \text{gph}F}} \frac{\langle y^*, (u-x) - (v-y) \rangle}{\|u-x\| + \|v-y\|} \le 0 \right\}$$

$$= \left\{ y^* \in X^* : \limsup_{\substack{(u,v) \to (x,y) \\ (u,v) \in \text{gph}F}} \frac{\langle y^*, (u-v) - (x-y) \rangle}{\|u-x\| + \|v-y\|} \le 0 \right\}. \qquad (3.1)$$

*In particular, let $f: \Delta \to X$ be a single valued and continuous mapping. Then,* for *$x \in \Delta$ and for any $y^* \in X^*$, we have*

$$\mathcal{F}\left(\widehat{D}^*f(x)\right)$$

$$= \left\{y^* \in X^*: \limsup_{u\to x} \frac{\langle y^*, u-x\rangle - \langle y^*, f(u)-f(x)\rangle}{\|u-x\|+\|f(u)-f(x)\|} \le 0\right\}$$

$$= \left\{y^* \in X^*: \limsup_{u\to x} \frac{\langle y^*,\ (u-x)-(f(u)-f(x))\rangle}{\|u-x\|+\|f(u)-f(x)\|} \le 0\right\}$$

$$= \left\{y^* \in X^*: \limsup_{u\to x} \frac{\langle y^*,\ (u-f(u))-(x-f(x))\rangle}{\|u-x\|+\|f(u)-f(x)\|} \le 0\right\}. \tag{3.2}$$

*Proof*. By definitions, the proof of this lemma is straight forward and it is omitted here. □

Fixed points of the Mordukhovich differential operator $\widehat{D}^*F(x, y)$ has the following basic properties.

**Proposition 3.2**. *Let* $\Delta$ *be a nonempty subset of a Banach space* $X$. *Let* $F: \Delta \rightrightarrows X$ *be a set valued mapping. Then, for any* $(x, y) \in \text{gph}F$, *we have*

(a) $\theta^* \in \mathcal{F}\left(\widehat{D}^*F(x, y)\right)$;

(b) $\mathcal{F}\left(\widehat{D}^*F(x, y)\right)$ *is a* $\|\cdot\|_*$*-closed and convex subset in* $X^*$.

*Proof*. Part (a) is clear. We prove (b). Let $(x, y) \in \text{gph}F$. By Part (a), $\mathcal{F}\left(\widehat{D}^*F(x, y)\right) \neq \emptyset$. For any $x^*, y^* \in \mathcal{F}\left(\widehat{D}^*F(x, y)\right)$ and for any $t \in [0, 1]$, by (3.1), we have

$$\limsup_{\substack{(u,v)\to(x,y)\\(u,v)\in\text{gph}F}} \frac{\langle tx^*+(1-t)y^*,\ (u-v)-(x-y)\rangle}{\|u-x\|+\|v-y\|}$$

$$\le t \limsup_{\substack{(u,v)\to(x,y)\\(u,v)\in\text{gph}F}} \frac{\langle x^*,\ (u-v)-(x-y)\rangle}{\|u-x\|+\|v-y\|} + (1-t) \limsup_{\substack{(u,v)\to(x,y)\\(u,v)\in\text{gph}F}} \frac{\langle y^*,\ (u-v)-(x-y)\rangle}{\|u-x\|+\|v-y\|}$$

$$\le 0.$$

By (3.1) again, this implies that $tx^* + (1-t)y^* \in \mathcal{F}\left(\widehat{D}^*F(x, y)\right)$. Hence, $\mathcal{F}\left(\widehat{D}^*F(x, y)\right)$ is convex. Next, we show that $\mathcal{F}\left(\widehat{D}^*F(x, y)\right)$ is closed in $(X^*, \|\cdot\|_*)$. Take an arbitrary $\|\cdot\|_*$-convergent sequence $\{y_n^*\} \subseteq \mathcal{F}\left(\widehat{D}^*F(x, y)\right)$ with limit $y^*$. For any $\varepsilon > 0$, there is $N \gg 1$ such that

$$\|y_n^* - y^*\|_* < \varepsilon, \text{ for any } n > N. \tag{3.3}$$

Fix an $m > N$. There is $\delta > 0$ with $\delta < 1$ such that

$$\frac{\langle y_m^*,\ (u-v)-(x-y)\rangle}{\|u-x\|+\|v-y\|} < \varepsilon, \text{ for any } (u, v) \in \text{gph}F \text{ with } \|u-x\| + \|v-y\| < \delta. \tag{3.4}$$

This implies that, for any $(u,v) \in \text{gph}F$ with $\|u-x\| + \|v-y\| < \delta$, by (3.3) and (3.4), we have

$$\frac{\langle y^*,\ (u-v)-(x-y)\rangle}{\|u-x\|+\|v-y\|}$$

$$= \frac{\langle y^*-y_m^*,\ (u-x)-(v-y)\rangle}{\|u-x\|+\|v-y\|} + \frac{\langle y_m^*,\ (u-v)-(x-y)\rangle}{\|u-x\|+\|v-y\|}$$

$$\leq \frac{\|y_m^*-y^*\|_*(\|u-x\|+\|v-y\|)}{\|u-x\|+\|v-y\|} + \varepsilon$$

$$< 2\varepsilon.$$

This implies

$$\limsup_{\substack{(u,v)\to(x,y)\\(u,v)\in\text{gph}F}} \frac{\langle y^*,\ (u-v)-(x-y)\rangle}{\|u-x\|+\|v-y\|} \leq 0.$$

By (3.1), this implies $y^* \in \mathcal{F}\left(\widehat{D}^*F(x,y)\right)$. □

**Proposition 3.3**. *Let Δ be a nonempty open subset of a reflexive Banach space X. Let $I_\Delta$ and $I_{X^*}$ be the identity mappings on Δ and $X^*$, respectively. Let $x_0 \in X$. Define a single valued continuous mapping f*: Δ → *X by*

$$f(x) = x_0 + x, \text{ for } x \in \Delta.$$

*Then, for any x ∈ Δ, we have*

(a) $\widehat{D}^*f(x)(y^*) = \{y^*\}$, *for any* $y^* \in X^*$;
(b) $\mathcal{F}\left(\widehat{D}^*f(x)\right) = X^*$;
(c) $\widehat{D}^*f(x) = I_{X^*}$.

*In particular, if $x_0 = \theta$, then f = $I_\Delta$. We have*

(d) $\widehat{D}^*I_\Delta(x)(y^*) = \{y^*\}$, *for any* $y^* \in X^*$;
(e) $\mathcal{F}\left(\widehat{D}^*I_\Delta(x)\right) = X^*$;
(f) $\widehat{D}^*I_\Delta(x) = I_{X^*}$.

*Proof*. It is clear that (b, c) follow from part (a) immediately. So, we only prove (a). As a matter of fact, part (a) can be proved by using the connection property between Mordukhovich derivatives and Fréchet derivatives for single valued mappings in Banach spaces (see Theorem 1.38 in [19]). However, we directly prove (a) here. Let *x* ∈ Δ. Since $f$ is a single valued continuous mapping, for any $y^* \in X^*$, we have

$$\limsup_{u\to x} \frac{\langle y^*,\ (u-x)-(u-x)\rangle}{\|u-x\|+\|u-x\|} = \limsup_{u\to x} \frac{\langle y^*,\ \theta\rangle}{\|u-x\|+\|u-x\|} = 0.$$

By (2.6), this implies

$$y^* \in \widehat{D}^* f(x)(y^*). \tag{3.5}$$

For any $x^* \in X^*$ with $x^* \neq y^*$, we calculate

$$\begin{aligned}
&\limsup_{\substack{u \to x \\ u \in \Delta}} \frac{\langle x^*, u - x\rangle - \langle y^*, f(u) - f(x)\rangle}{\|u - x\| + \|f(u) - f(x)\|} \\
&= \limsup_{\substack{u \to x \\ u \in \Delta}} \frac{\langle x^*, u - x\rangle - \langle y^*, u - x\rangle}{\|u - x\| + \|u - x\|} \\
&= \limsup_{\substack{u \to x \\ u \in \Delta}} \frac{\langle x^* - y^*, u - x\rangle}{\|u - x\| + \|u - x\|}.
\end{aligned} \tag{3.6}$$

Since $X$ is a reflexive Banach space and $x^* - y^* \neq \theta^*$, there is $j^*(x^* - y^*) \in J^*(x^* - y^*) \subseteq X$ satisfying

$$\langle x^* - y^*,\ j^*(x^* - y^*)\rangle = \|j^*(x^* - y^*)\|^2 = \|x^* - y^*\|_*^2 > 0. \tag{3.7}$$

In the limit (3.5), we take a special direction $u_t = x + tj^*(x^* - y^*)$, for $t > 0$. Since $\Delta$ is open and $x \in \Delta$, there is $a > 0$ such that

$$u_t = x + tj^*(x^* - y^*) \in \Delta, \text{ for all } t \text{ satisfying } a > t > 0.$$

Then, by (3.7), we have

$$\begin{aligned}
&\limsup_{\substack{u \to x \\ u \in \Delta}} \frac{\langle x^* - y^*, u - x\rangle}{\|u - x\| + \|u - x\|} \\
&\geq \limsup_{\substack{t \downarrow 0 \\ t < a}} \frac{\langle x^* - y^*, x + tj^*(x^* - y^*) - x\rangle}{\|x + tj^*(x^* - y^*) - x\| + \|x + tj^*(x^* - y^*) - x\|} \\
&= \limsup_{\substack{t \downarrow 0 \\ t < a}} \frac{t\|j^*(x^* - y^*)\|^2}{2t\|j^*(x^* - y^*)\|} \\
&= \frac{\|j^*(x^* - y^*)\|}{2} > 0.
\end{aligned}$$

This implies that

$$x^* \notin \widehat{D}^* f(x)(y^*), \text{ for any } x^* \in X^* \text{ with } x^* \neq y^*. \tag{3.8}$$

Then, (a) is proved by (3.5) and (3.8). □

By Proposition 3.3, we have the following observations.

**Remark 3.4**. *Let $\Delta$ be a nonempty open subset of a reflexive Banach space $X$. Let $I_\Delta$ and $I_{X^*}$ be*

*the identity mappings on* $\Delta$ *and* $X^*$*, respectively. Let* $f$ *be a single valued continuous mapping* $f\colon \Delta \to X$*. Then,*

(i) $f = I_\Delta \quad \Longrightarrow \quad \widehat{D}^* f(x) = I_{X^*}$, *for any* $x \in \Delta$;

(ii) $\widehat{D}^* f(x) = I_{X^*}$, *for any* $x \in \Delta \quad \not\Rightarrow \quad f = I_\Delta$.

**Proposition 3.5**. *Let* $\Delta$ *be a nonempty open subset of a reflexive Banach space* $X$*. Let* $x_0 \in X$ *with* $x_0 \neq \theta$ *and* $\lambda \neq 1$*. Define a single-valued continuous mapping* $f\colon \Delta \to X$ *by*

$$f(x) = x_0 + \lambda x, \text{ for } x \in \Delta.$$

*Then,* $\mathcal{F}\left(\widehat{D}^* f(x)\right) = \{\theta^*\}$*, for any* $x \in \Delta$.

*Proof.* The proof of this lemma is similar to the proof of part (a) of Proposition 3.3. For any $x^* \in X^*$ with $x^* \neq \theta^*$, we calculate

$$\limsup_{\substack{u\to x\\ u\in\Delta}} \frac{\langle x^*, u-x\rangle - \langle x^*, f(u)-f(x)\rangle}{\|u-x\|+\|f(u)-f(x)\|}$$

$$= \limsup_{\substack{u\to x\\ u\in\Delta}} \frac{\langle x^*, u-x\rangle - \langle x^*, \lambda(u-x)\rangle}{\|u-x\|+|\lambda|\|u-x\|}$$

$$= \limsup_{\substack{u\to x\\ u\in\Delta}} \frac{\langle x^*, (1-\lambda)(u-x)\rangle}{(1+|\lambda|)\|u-x\|}. \tag{3.9}$$

Since $X$ is a reflexive Banach space and $x^* \neq \theta^*$, there is $j^*(x^*) \in J^*(x^*) \subseteq X$ satisfying

$$\langle x^*,\ j^*(x^*)\rangle = \|j^*(x^*)\|^2 = \|x^*\|_*^2 > 0. \tag{3.10}$$

Case 1. $\lambda < 1$. In this case, in the limit (3.9), we take a special direction $u_t = x + tj^*(x^*)$, for $t > 0$. Since $\Delta$ is open and $x \in \Delta$, there is $a > 0$ such that

$$u_t = x + tj^*(x^* - y^*) \in \Delta, \text{ for all } t \text{ satisfying } a > t > 0.$$

Then, by (3.10), we have

$$\limsup_{\substack{u\to x\\ u\in\Delta}} \frac{\langle x^*, (1-\lambda)(u-x)\rangle}{(1+|\lambda|)\|u-x\|}$$

$$\geq \limsup_{\substack{t\downarrow 0\\ t<a}} \frac{(1-\lambda)(\langle x^*,\ x+tj^*(x^*)-x\rangle}{(1+|\lambda|)t\|j^*(x^*)\|}$$

$$= \limsup_{\substack{t\downarrow 0\\ t<a}} \frac{(1-\lambda)t\|j^*(x^*)\|^2}{(1+|\lambda|)t\|j^*(x^*-y^*)\|}$$

$$= \frac{(1-\lambda)\|j^*(x^*-y^*)\|}{1+|\lambda|} > 0. \tag{3.11}$$

Case 2. λ > 1. In this case, in the limit (3.9), we take a special direction $u_t = x - tj^*(x^*)$, for $t$ > 0. Since Δ is open and $x \in \Delta$, there is $a$ > 0 such that

$$u_t = x - tj^*(x^* - y^*) \in \Delta, \text{ for all } t \text{ satisfying } a > t > 0.$$

Then, by (3.10), we have

$$\limsup_{\substack{u\to x\\ u\in\Delta}} \frac{\langle x^*, (1-\lambda)(u-x)\rangle}{(1+|\lambda|)\|u-x\|}$$

$$\geq \limsup_{\substack{t\downarrow 0\\ t<a}} \frac{(1-\lambda)(\langle x^*,\ x-tj^*(x^*)-x\rangle}{(1+|\lambda|)t\|j^*(x^*)\|}$$

$$= \limsup_{\substack{t\downarrow 0\\ t<a}} \frac{-(1-\lambda)t\|j^*(x^*)\|^2}{(1+|\lambda|)t\|j^*(x^*-y^*)\|}$$

$$= \frac{-(1-\lambda)\|j^*(x^*-y^*)\|}{1+|\lambda|} > 0. \tag{3.12}$$

(3.11) and (3.12) together imply that

$$x^* \notin \widehat{D}^* I_\Delta(x)(x^*), \text{ for any } x^* \in X^* \text{ with } x^* \neq \theta^*. \tag{3.13}$$

Then, this proposition is proved by (3.13) and Lemma 3.1. □

## 4. Some fixed-point properties of the Mordukhovich differential operator with respect to the standard metric projection

In this section, we consider the standard metric projection, which is a special mapping in Banach spaces. This operator has been widely applied in optimization theory, approximation theory, fixed point theory, and so forth (see [1, 5, 8, 11, 20−21, 24]). We use some results about Mordukhovich derivatives of the metric projection obtained in [13−15] to investigate fixed-point properties of the Mordukhovich differential operator with respect to the standard metric projection.

### 4.1. The metric projection onto closed balls in Banach spaces

In this subsection, let $X$ be a uniformly convex and uniformly smooth Banach space with closed unit ball $\mathbb{B}$. For any $r$ > 0, $r\mathbb{B}$ is the closed ball in $X$ with center θ and radius $r$. The metric projection $P_{r\mathbb{B}}: X \to r\mathbb{B}$ has the following analytic representations (see [14]).

$$P_{r\mathbb{B}}(x) = \begin{cases} x, & \text{for any } x \in r\mathbb{B}, \\ \frac{r}{\|x\|}x, & \text{for any } x \notin r\mathbb{B}. \end{cases}$$

**Theorem 3.1 (partial) in [14].** *Let X be a uniformly convex and uniformly smooth Banach space. For any r > 0, the Mordukhovich derivatives of the metric projection* $P_{r\mathbb{B}}: X \to r\mathbb{B}$ *has the following properties.*

(i) *For every* $\bar{x} \in r\mathbb{B}^{o}$, *we have*

$$\widehat{D}^* P_{r\mathbb{B}}(\bar{x})(y^*) = \{y^*\}, \text{ for every } y^* \in X^*.$$

(ii) *For every* $\bar{x} \in X \backslash r\mathbb{B}$, *we have*

$$\widehat{D}^* P_{r\mathbb{B}}(\bar{x})(y^*) = \left\{\frac{r}{\|\bar{x}\|}\left(y^* - \frac{\langle y^*, \bar{x}\rangle}{\|\bar{x}\|^2} J(\bar{x})\right)\right\}, \text{ for every } y^* \in X^*.$$

*In particular, we have*

(a) $\widehat{D}^* P_{r\mathbb{B}}(\bar{x})(y^*) = \left\{\frac{r}{\|\bar{x}\|} y^*\right\}$, *if* $y^* \perp \bar{x}$;

(b) $\widehat{D}^* P_{r\mathbb{B}}(\bar{x})(J(\bar{x})) = \{\theta^*\}$.

**Theorem 4.1.** *Let X be a uniformly convex and uniformly smooth Banach space. For any r* > 0, *we have*

(i) $$\widehat{D}^* P_{r\mathbb{B}}(x) = I_{X^*} \quad \text{and} \quad \mathcal{F}\left(\widehat{D}^* P_{r\mathbb{B}}(x)\right) = X^*, \text{ for every } x \in r\mathbb{B}^{o}.$$

(ii) $$\mathcal{F}\left(\widehat{D}^* P_{r\mathbb{B}}(x)\right) = \{\theta^*\}, \text{ for any } x \in X \backslash r\mathbb{B}.$$

*Proof*. Proof of (i). Since $P_{r\mathbb{B}} = I_{r\mathbb{B}}$ on $r\mathbb{B}^{o}$ and $r\mathbb{B}^{o}$ is open, then, part (i) of this theorem follows from Proposition 3.3 immediately.

Proof of (ii). By Theorem 3.1 in [14], for every $x \in X \backslash r\mathbb{B}$, we have

$$\widehat{D}^* P_{r\mathbb{B}}(x)(y^*) = \left\{\frac{r}{\|x\|}\left(y^* - \frac{\langle y^*, x\rangle}{\|x\|^2} J(x)\right)\right\}, \text{ for every } y^* \in X^*.$$

Solve for $y^*$ in the following equation

$$\frac{r}{\|x\|}\left(y^* - \frac{\langle y^*, x\rangle}{\|x\|^2} J(x)\right) = y^*. \tag{4.1}$$

Taking the inner product with *x* on both sides of the above equation, we obtain

$$\frac{r}{\|x\|}(\langle y^*, x\rangle - \langle y^*, x\rangle) = \langle y^*, x\rangle.$$

This implies $\langle y^*, x\rangle = 0$. Substituting it into (4.1), by $\frac{r}{\|x\|} < 1$, we have $y^* = \theta^*$. □

### 4.2 The metric projection onto the positive cone in $L_p(S)$

In this subsection, we first recall some notations and definitions about the real uniformly convex and uniformly smooth Banach space $(L_p(S), \|\cdot\|_p)$, which are used in [14]. Let $(S, \mathcal{A}, \mu)$ be a positive and complete measure space. The real uniformly convex and uniformly smooth Banach space $(L_p(S), \|\cdot\|_p)$ has dual space $(L_q(S), \|\cdot\|_q)$, in which $p$ and $q$ satisfy $1 < p, q < \infty$ and $\frac{1}{p} + \frac{1}{q} = 1$. $L_p(S)$ and $L_q(S)$ share the same origin $\theta = \theta^*$(see [5] for more details). For the sake of distinction between $L_p(S)$ and its dual space $L_q(S)$, we use English letters $f, g, h, \ldots$ for the elements in $L_p(S)$, and we use Greek letters $\varphi, \psi, \xi, \ldots$ for the elements in the dual space $L_q(S)$. The positive cone $K_p$ and $K_q$ of $L_p(S)$ and $L_q(S)$ are respectively defined by

$$K_p = \{f \in L_p(S) \colon f(s) \geq 0, \text{ for } \mu\text{-almost all } s \in S\},$$

$$K_q = \{\varphi \in L_q(S) \colon \varphi(s) \geq 0, \text{ for } \mu\text{-almost all } s \in S\}.$$

Let $\preccurlyeq_p$ and $\preccurlyeq_q$ be the two partial orders on $L_p(S)$ and $L_q(S)$ induced by $K_p$ and $K_q$, respectively. For any $f, g \in L_p(S)$ with $f \preccurlyeq_p g$, and for any $\xi, \psi \in L_q(S)$ with $\xi \preccurlyeq_q \psi$, we write

$$[f, g]_{\preccurlyeq_p} = \{h \in L_p(S) \colon f \preccurlyeq_p h \preccurlyeq_p g\},$$

$$[\xi, \psi]_{\preccurlyeq_q} = \{\varphi \in L_q(S) \colon \xi \preccurlyeq_q \varphi \preccurlyeq_q \psi\}.$$

The normalized duality mapping $J \colon L_p(S) \to L_q(S)$ has the following representation, for any given $f \in L_p(S)$ with $f \neq \theta$,

$$(Jf)(s) = \frac{|f(s)|^{p-1}\mathrm{sign}(f(s))}{\|f\|_p^{p-2}} = \frac{|f(s)|^{p-2}f(s)}{\|f\|_p^{p-2}}, \text{ for all } s \in S.$$

**Lemma 5.2 in [14]**. *The metric projection $P_{K_p} \colon L_p(S) \to K_p$ has the following representations.*

$$(P_{K_p}f)(s) = \begin{cases} f(s), & \text{if } f(s) > 0, \\ 0, & \text{if } f(s) \leq 0, \end{cases} \quad \text{for any } f \in L_p(S). \tag{4.2}$$

**Theorem 5.2 in [14].** *The metric projection $P_{K_p} \colon L_p(S) \to K_p$ has the following Mordukhovich derivatives. For any $f \in L_p(S)$,*

(a) $\widehat{D}^* P_{K_p}(f)(\theta^*) = \{\theta^*\}$;

(b) *For any $\varphi \in L_q(S)$,*

$$\theta^* \in \widehat{D}^* P_{K_p}(f)(\varphi)$$

$$\Leftrightarrow \quad \mu(\{s \in S \colon \varphi(s) \neq 0 \ \text{ and } f(s) > 0\} \cup \{s \in S \colon \varphi(s) < 0 \ \text{ and } f(s) \leq 0\}) = 0.$$

*This implies that,*

($b_1$) *For any $f \in -K_p$, we have*

$$\theta^* \in \widehat{D}^* P_{K_p}(f)(\varphi), \text{ for any } \varphi \in K_q;$$

($b_2$) *For any $f \in K_p \backslash \{\theta\}$,*

$$\theta^* \notin \widehat{D}^* P_{K_p}(f)(J(f));$$

(c) $J(f) \in \widehat{D}^* P_{K_p}(f)(J(f))$, *for any* $f \in K_p$;

(d) *For any given* $\psi \in K_q$, *we have*

$$\widehat{D}^* P_{K_p}(\theta)(\psi) = [\theta^*, \psi]_{\preccurlyeq_q}.$$

**Theorem 4.2.** *The Mordukhovich differential operator* $\widehat{D}^* P_{K_p}(\cdot)$ *of the metric projection* $P_{K_p}: L_p(S) \to K_p$ *has the following fixed-point properties.*

(i) $J(f) \in \mathcal{F}\left(\widehat{D}^* P_{K_p}(f)\right)$, *for any* $f \in K_p$;

(ii) $\psi \in \mathcal{F}\left(\widehat{D}^* P_{K_p}(\theta)\right)$, *for any* $\in K_q$.

*Proof.* Parts (i) and (ii) of this theorem follow from parts (c) and (d) of Theorem 5.2 in [14], respectively. □

### 4.3 The metric projection in $l_2$

In this subsection, we consider the real Hilbert space $l_2$ with norm $\|\cdot\|$, inner product $\langle\cdot, \cdot\rangle$ and the origin $\theta$. Let $\mathbb{N}$ denote the set of all positive integers. We recall some definitions and notations in the real Hilbert space $l_2$ used in [13].

$\mathbb{K} = \{x = (x_1, x_2, \ldots) \in l_2: x_i \geq 0, \text{ for all } i \in \mathbb{N}\}$;

$\mathbb{K}^+ = \{x = (x_1, x_2, \ldots) \in \mathbb{K}: x_i > 0, \text{ for all } i \in \mathbb{N}\}$;

$\mathbb{K}^- = \{x = (x_1, x_2, \ldots) \in l_2: x_i < 0, \text{ for all } i \in \mathbb{N}\}$;

$\widehat{\mathbb{K}} = \{x = (x_1, x_2, \ldots) \in l_2: |x_i| > 0, \text{ for all } i \in \mathbb{N}, \text{ and there are } j, k \in \mathbb{N} \text{ with } x_j x_k < 0\}$.

Let $N$ be a nonempty subset of $\mathbb{N}$ with complement $\overline{N}$. We define some subsets in $l_2$ with respect to $N$.

$\mathbb{R}^N = \{x = (x_1, x_2, \ldots) \in l_2: x_i = 0, \text{ for all } i \in \overline{N}\}$;

$\mathbb{K}_N = \{x = (x_1, x_2, \ldots) \in l_2: x_i \geq 0, \text{ for all } i \in N\}$;

$\partial\mathbb{K}_N = \{x = (x_1, x_2, \ldots) \in \mathbb{K}_N: x_i = 0, \text{ for some } i \in N\}$;

$\mathbb{Z}_N = \{x = (x_1, x_2, \ldots) \in l_2: x_i > 0, \text{ for all } i \in N, \text{ and } x_i = 0, \text{ for all } i \in \overline{N}\}$.

$\mathbb{R}^N$ is a closed subspace of $l_2$. $\mathbb{K}$ and $\mathbb{K}_N$ all are pointed closed and convex cones in $l_2$. We define an ordering relation $\preccurlyeq_N$ on $l_2$, for any $y, z \in l_2$, by

$$z \preccurlyeq_N y \quad \Leftrightarrow \quad z_i \leq y_i, \text{ for } i \in N \text{ and } z_i = y_i, \text{ for } i \in \overline{N}.$$

**Lemma 5.2 in [13]**. *Let $\mathbb{K}$ be the positive cone of $l_2$. $P_{\mathbb{K}}$ has the following properties.*

(a) *For any $x \in l_2$, $P_{\mathbb{K}}(x)$ is represented as follows*

$$P_{\mathbb{K}}(x)_i = \begin{cases} x_i, & \text{if } x_i > 0, \\ 0, & \text{if } x_i \le 0, \end{cases} \quad \text{for } i \in \mathbb{N}.$$

(b) *$P_{\mathbb{K}}$ is positive homogeneous. For any $x \in l_2$,*

$$P_{\mathbb{K}}(\lambda x) = \lambda P_{\mathbb{K}}(x), \text{ for any } \lambda \ge 0.$$

**Theorem 5.1 in [13]**. *Let $M$ be a nonempty finite subset of $\mathbb{N}$ with complement $\overline{M}$. The Mordukhovich derivatives of $P_{\mathbb{K}}$ have the following representations.*

(i) *For any $\bar{x} \in l_2$, we have*

$$\widehat{D}^* P_{\mathbb{K}}(\bar{x})(\theta) = \{\theta\};$$

(ii) *Let $\bar{x} \in \mathbb{Z}_M$. For any $y \in l_2$,*

$$y \in \mathbb{K}_{\overline{M}} \iff y \in \widehat{D}^* P_{\mathbb{K}}(\bar{x})(y).$$

(iii) *Let $\bar{x} \in \mathbb{Z}_M$. For any $y \in \mathbb{K}_{\overline{M}}$, we have*

$$\widehat{D}^* P_{\mathbb{K}}(\bar{x})(y) = \{z \in \mathbb{K}_{\overline{M}} : z \preccurlyeq_{\overline{M}} y\}.$$

*In particular,*

$$y \in \partial \mathbb{K}_{\overline{M}} \implies \widehat{D}^* P_{\mathbb{K}}(\bar{x})(y) = \{y\}.$$

**Theorem 4.3**. *Let $M$ be a nonempty finite subset of $\mathbb{N}$ with complement $\overline{M}$. Then, for any $\bar{x} \in \mathbb{Z}_M$, we have*

$$\mathcal{F}\left(\widehat{D}^* P_{\mathbb{K}}(\bar{x})\right) = \mathbb{K}_{\overline{M}}.$$

*Proof.* This theorem follows parts (ii) and (iii) of Theorem 5.1 in [13]. □

### 4.4. The metric projection operator in the real Banach space $l_1$

In this subsection, we consider the real Banach space $(l_1, \|\cdot\|_1)$ with dual space $(l_\infty, \|\cdot\|_\infty)$. It is well-known that $l_1$ is neither uniformly convex, nor uniformly smooth. For any $r > 0$, let $r\mathbb{B}$ be the closed ball in $l_1$ centered at θ and with radius $r$. In [15], it is proved that the metric projection $P_{r\mathbb{B}}$ is indeed a set valued mapping (not a single valued mapping). It causes that both the Gâteaux directional differentiability and Fréchet differentiability of $P_{r\mathbb{B}}$ are not defined in $l_1$. In [15], the Mordukhovich derivative of $P_{r\mathbb{B}}$ in $l_1$ was studied. In this subsection, we use the results obtained in [15] to study the fixed-point properties of the Mordukhovich differential operator of the metric projection $P_{r\mathbb{B}}$ in $l_1$, which is a set valued mapping.

**Proposition 3.2 in [15]**. *Let $r > 0$. For any $x \in l_1$ with $\|x\|_1 > r$, the projection $P_{r\mathbb{B}}$ satisfies*

(a) $\frac{r}{\|x\|_1}x \in P_{r\mathbb{B}}(x)$;

(b) *For any* $j(x) \in J(x)$, *we have*

$$\langle j(x),\ \frac{r}{\|x\|_1}x - y\rangle \geq 0, \textit{ for all } y \in r\mathbb{B}.$$

**Theorem 3.4 in [15].** *For any* $r > 0$, *the Mordukhovich derivatives of the set valued metric projection* $P_{r\mathbb{B}}$: $l_1 \rightrightarrows r\mathbb{B}$ *has the following properties.*

(i) *For every* $x \in r\mathbb{B}^{\mathrm{o}}$ *with* $\|x\|_1 < r$, *we have*

$$\widehat{D}^*P_{r\mathbb{B}}(x,x)(\varphi) = \{\varphi\}, \textit{ for every } \varphi \in l_\infty.$$

(ii) *For every* $x \in K_1^+$ *with* $\|x\|_1 > r$, *we have*

(a) $\widehat{D}^*P_{r\mathbb{B}}\left(x, \frac{r}{\|x\|_1}x\right)(\theta^*) = \{\theta^*\}$;

(b) $\widehat{D}^*P_{r\mathbb{B}}\left(x, \frac{r}{\|x\|_1}x\right)(\beta_d) = \emptyset$, *for any* $d > 0$;

(c) $\widehat{D}^*P_{r\mathbb{B}}\left(x, \frac{r}{\|x\|_1}x\right)(J(x)) = \emptyset$.

**Theorem 4.4.** *Let* $r > 0$. *The Mordukhovich differential operator* $\widehat{D}^*P_{r\mathbb{B}}(\cdot)$ *of the metric projection* $P_{r\mathbb{B}}$:$l_1 \to r\mathbb{B}$ *has the following fixed-point properties.*

(i) *For every* $x \in r\mathbb{B}^{\mathrm{o}}$, *we have*

$$\widehat{D}^*P_{r\mathbb{B}}\left(x, \frac{r}{\|x\|_1}x\right) = I_{l_\infty} \quad \textit{and} \quad \mathcal{F}\left(\widehat{D}^*P_{r\mathbb{B}}\left(x, \frac{r}{\|x\|_1}x\right)\right) = l_\infty.$$

(ii) *For any* $x \in l_1 \backslash r\mathbb{B}$, *we have*

$$\mathcal{F}\left(\widehat{D}^*P_{r\mathbb{B}}\left(x, \frac{r}{\|x\|_1}x\right)\right) = \{\theta^*\}.$$

*Proof*. Proof of (i). Since $P_{r\mathbb{B}} = I_{r\mathbb{B}}$ on $r\mathbb{B}^{\mathrm{o}}$ and $r\mathbb{B}^{\mathrm{o}}$ is open, then, part (i) of this theorem follows from Proposition 3.3 immediately.

Proof of (ii). Let $x \in l_1 \backslash r\mathbb{B}$ with $x = (x_1, x_1, \ldots)$, by Proposition 3.3, we have

$$\theta^* \in \mathcal{F}\left(\widehat{D}^*P_{r\mathbb{B}}\left(x, \frac{r}{\|x\|_1}x\right)\right). \tag{4.3}$$

Next, we prove

$$\varphi \notin \mathcal{F}\left(\widehat{D}^*P_{r\mathbb{B}}\left(x, \frac{r}{\|x\|_1}x\right)\right), \text{ for any } \varphi \in l_\infty \text{ with } \varphi \neq \theta^*. \tag{4.4}$$

For any $\varphi \in l_\infty$ with $\varphi \neq \theta^*$, there is a positive number $n$ such that the $n$th coordinate of $\varphi$, denoted by $\varphi_n$, is not zero. Define $y \in l_1$ satisfying that the $n$th coordinate of $y$ equals $\varphi_n$ and all other coordinates of $y$ are zero. Since $l_1 \backslash r\mathbb{B}$ is open and $x \in l_1 \backslash r\mathbb{B}$, there is $a > 0$ such that

$$u_t = x + ty \in l_1 \backslash r\mathbb{B}, \text{ for all } t \text{ satisfying } a > t > 0.$$

By Proposition 3.2 in [15], we have

$$\frac{r}{\|u_t\|_1} u_t \in P_{r\mathbb{B}}(u_t), \text{ for } u_t = x + ty \in l_1 \backslash r\mathbb{B}, \text{ for all } t \text{ satisfying } a > t > 0.$$

Then, we have

$$\limsup_{\substack{(u,v)\to\left(x,\frac{r}{\|x\|_1}x\right)\\ v\in P_{r\mathbb{B}}(u)}} \frac{\langle\varphi, u-x\rangle - \langle\varphi,\ v-\frac{r}{\|x\|_1}x\rangle}{\|u-x\|_1 + \left\|v-\frac{r}{\|x\|_1}x\right\|_1}$$

$$\geq \limsup_{\left(u,\frac{r}{\|u\|_1}u\right)\to\left(x,\frac{r}{\|x\|_1}x\right)} \frac{\langle\varphi, u-x\rangle - \langle\varphi,\ \frac{r}{\|u\|_1}u - \frac{r}{\|x\|_1}x\rangle}{\|u-x\|_1 + \left\|\frac{r}{\|u\|_1}u - \frac{r}{\|x\|_1}x\right\|_1}$$

$$\geq \limsup_{\substack{t\downarrow 0\\ t<a}} \frac{\langle\varphi, u_t-x\rangle - \langle\varphi,\ \frac{r}{\|u_t\|_1}u_t - \frac{r}{\|x\|_1}x\rangle}{\|u_t-x\|_1 + \left\|\frac{r}{\|u_t\|_1}u_t - \frac{r}{\|x\|_1}x\right\|_1}$$

$$= \limsup_{\substack{t\downarrow 0\\ t<a}} \frac{\langle\varphi, ty\rangle - \langle\varphi,\ \frac{r}{\|u_t\|_1}u_t - \frac{r}{\|x\|_1}x\rangle}{\|ty\|_1 + \left\|\frac{r}{\|u_t\|_1}u_t - \frac{r}{\|x\|_1}x\right\|_1}$$

$$= \limsup_{\substack{t\downarrow 0\\ t<a}} \frac{t\varphi_n^2 - \langle\varphi,\ \frac{r}{\|u_t\|_1}(x+ty) - \frac{r}{\|x\|_1}x\rangle}{t|\varphi_n| + \left\|\frac{r}{\|u_t\|_1}u_t - \frac{r}{\|x\|_1}x\right\|_1}$$

$$= \limsup_{\substack{t\downarrow 0\\ t<a}} \frac{t\varphi_n^2 - \langle\varphi,\ \frac{r}{\|u_t\|_1}ty + \frac{r}{\|u_t\|_1}x - \frac{r}{\|x\|_1}x\rangle}{t|\varphi_n| + \left\|\frac{r}{\|u_t\|_1}u_t - \frac{r}{\|x\|_1}x\right\|_1}$$

$$= \limsup_{\substack{t\downarrow 0\\ t<a}} \frac{t\varphi_n^2 - \langle\varphi,\ \frac{r}{\|u_t\|_1}ty + \frac{r}{\|u_t\|_1}x - \frac{r}{\|x\|_1}x\rangle}{t|\varphi_n| + \left\|\frac{r}{\|u_t\|_1}ty + \frac{r}{\|u_t\|_1}x - \frac{r}{\|x\|_1}x\right\|_1}$$

$$\geq \limsup_{\substack{t\downarrow 0\\ t<a}} \frac{t\varphi_n^2 - \langle\varphi,\ \frac{r}{\|u_t\|_1}ty + \frac{r}{\|u_t\|_1}x - \frac{r}{\|x\|_1}x\rangle}{t|\varphi_n| + \left\|\frac{r}{\|u_t\|_1}ty\right\|_1 + \left\|\frac{r}{\|u_t\|_1}x - \frac{r}{\|x\|_1}x\right\|_1}$$

$$= \limsup_{\substack{t\downarrow 0\\ t<a}} \frac{t\varphi_n^2\left(1-\frac{r}{\|u_t\|_1}\right)-\left\langle \varphi, \left(\frac{r}{\|u_t\|_1}-\frac{r}{\|x\|_1}\right)x\right\rangle}{t|\varphi_n|+t|\varphi_n|\frac{r}{\|u_t\|_1}+\left\|\frac{r}{\|u_t\|_1}x-\frac{r}{\|x\|_1}x\right\|_1}$$

$$= \limsup_{\substack{t\downarrow 0\\ t<a}} \frac{t\varphi_n^2\left(1-\frac{r}{\|u_t\|_1}\right)-\langle \varphi, x\rangle\left(\frac{r}{\|u_t\|_1}-\frac{r}{\|x\|_1}\right)}{t|\varphi_n|\left(1+\frac{r}{\|u_t\|_1}\right)+r\|x\|_1\left|\frac{1}{\|u_t\|_1}-\frac{1}{\|x\|_1}\right|}. \qquad (4.5)$$

The proof of (4.4) is divided into the following three cases.

Case 1. Suppose $\langle \varphi,\ x\rangle = 0$. Then, the limit (4.5) is calculated as follows.

$$\limsup_{\substack{t\downarrow 0\\ t<a}} \frac{t\varphi_n^2\left(1-\frac{r}{\|u_t\|_1}\right)-\langle \varphi, x\rangle\left(\frac{r}{\|u_t\|_1}-\frac{r}{\|x\|_1}\right)}{t|\varphi_n|\left(1+\frac{r}{\|u_t\|_1}\right)+r\|x\|_1\left|\frac{1}{\|u_t\|_1}-\frac{1}{\|x\|_1}\right|}$$

$$= \limsup_{\substack{t\downarrow 0\\ t<a}} \frac{t\varphi_n^2\left(1-\frac{r}{\|u_t\|_1}\right)}{t|\varphi_n|\left(1+\frac{r}{\|u_t\|_1}\right)+r\|x\|_1\left|\frac{|x_n|-|x_n+t\varphi_n|}{\|x\|_1\|u_t\|_1}\right|}$$

$$= \limsup_{\substack{t\downarrow 0\\ t<a}} \frac{t\varphi_n^2\left(1-\frac{r}{\|u_t\|_1}\right)}{t|\varphi_n|\left(1+\frac{r}{\|u_t\|_1}\right)+r\frac{||x_n|-|x_n+t\varphi_n||}{\|u_t\|_1}}$$

$$= \limsup_{\substack{t\downarrow 0\\ t<a}} \frac{t\varphi_n^2\left(1-\frac{r}{\|u_t\|_1}\right)}{t|\varphi_n|\left(1+\frac{r}{\|u_t\|_1}\right)+r\frac{|-2t\varphi_n-t^2\varphi_n^2|}{\|u_t\|_1(|x_n|+|x_n+t\varphi_n|)}}$$

$$= \limsup_{\substack{t\downarrow 0\\ t<a}} \frac{|\varphi_n|\left(1-\frac{r}{\|u_t\|_1}\right)}{\left(1+\frac{r}{\|u_t\|_1}\right)+r\frac{|2+t\varphi_n|}{\|u_t\|_1(|x_n|+|x_n+t\varphi_n|)}}$$

$$= \frac{|\varphi_n|\left(1-\frac{r}{\|x\|_1}\right)}{\left(1+\frac{r}{\|x\|_1}\right)+\frac{r}{|x_n|\,\|x\|_1}}$$

$$> 0. \qquad (4.6)$$

Here, we assumed $x_n \neq 0$. In case if $x_n = 0$, then

$$\limsup_{\substack{t\downarrow 0\\ t<a}} \frac{t\varphi_n^2\left(1-\frac{r}{\|u_t\|_1}\right)}{t|\varphi_n|\left(1+\frac{r}{\|u_t\|_1}\right)+r\frac{||x_n|-|x_n+t\varphi_n||}{\|u_t\|_1}} = \frac{|\varphi_n|\left(1-\frac{r}{\|x\|_1}\right)}{\left(1+\frac{r}{\|x\|_1}\right)+\frac{r}{\|x\|_1}} > 0.$$

Suppose $\langle \varphi,\ x\rangle \neq 0$. The limit (4.5) is calculated as below.

$$\limsup_{\substack{t\downarrow 0\\ t<a}} \frac{t\varphi_n^2\left(1-\frac{r}{\|u_t\|_1}\right)-\langle\varphi,x\rangle\left(\frac{r}{\|u_t\|_1}-\frac{r}{\|x\|_1}\right)}{t|\varphi_n|\left(1+\frac{r}{\|u_t\|_1}\right)+r\|x\|_1\left|\frac{1}{\|u_t\|_1}-\frac{1}{\|x\|_1}\right|}$$

$$=\limsup_{\substack{t\downarrow 0\\ t<a}} \frac{t\varphi_n^2\left(1-\frac{r}{\|u_t\|_1}\right)-r\langle\varphi,x\rangle\frac{|x_n|-|x_n+t\varphi_n|}{\|x\|_1\|u_t\|_1}}{t|\varphi_n|\left(1+\frac{r}{\|u_t\|_1}\right)+r\|x\|_1\left|\frac{|x_n|-|x_n+t\varphi_n|}{\|x\|_1\|u_t\|_1}\right|}$$

$$=\limsup_{\substack{t\downarrow 0\\ t<a}} \frac{t\varphi_n^2\left(1-\frac{r}{\|u_t\|_1}\right)-r\langle\varphi,x\rangle\frac{-2tx_n\varphi_n-t^2\varphi_n^2}{\|x\|_1\|u_t\|_1(|x_n|+|x_n+t\varphi_n|)}}{t|\varphi_n|\left(1+\frac{r}{\|u_t\|_1}\right)+r\|x\|_1\left|\frac{-2tx_n\varphi_n-t^2\varphi_n^2}{\|x\|_1\|u_t\|_1(|x_n|+|x_n+t\varphi_n|)}\right|}. \tag{4.7}$$

Case 2. If $\langle\varphi,\ x\rangle > 0$, we can choose *n* such that $\varphi_n \neq 0$ and $x_n\varphi_n > 0$. Then, by $x_n\varphi_n > 0$, we calculate (4.7).

$$\limsup_{\substack{t\downarrow 0\\ t<a}} \frac{t\varphi_n^2\left(1-\frac{r}{\|u_t\|_1}\right)-r\langle\varphi,x\rangle\frac{-2tx_n\varphi_n-t^2\varphi_n^2}{\|x\|_1\|u_t\|_1(|x_n|+|x_n+t\varphi_n|)}}{t|\varphi_n|\left(1+\frac{r}{\|u_t\|_1}\right)+r\|x\|_1\left|\frac{-2tx_n\varphi_n-t^2\varphi_n^2}{\|x\|_1\|u_t\|_1(|x_n|+|x_n+t\varphi_n|)}\right|}$$

$$=\limsup_{\substack{t\downarrow 0\\ t<a}} \frac{t\varphi_n^2\left(1-\frac{r}{\|u_t\|_1}\right)+\frac{2tr\langle\varphi,x\rangle x_n\varphi_n}{\|x\|_1\|u_t\|_1(|x_n|+|x_n+t\varphi_n|)}+\frac{t^2r\langle\varphi,x\rangle\varphi_n^2}{\|x\|_1\|u_t\|_1(|x_n|+|x_n+t\varphi_n|)}}{t|\varphi_n|\left(1+\frac{r}{\|u_t\|_1}\right)+r\|x\|_1t|\varphi_n|\left|\frac{2x_n+t|\varphi_n|}{\|x\|_1\|u_t\|_1(|x_n|+|x_n+t\varphi_n|)}\right|}$$

$$=\limsup_{\substack{t\downarrow 0\\ t<a}} \frac{t\varphi_n^2\left(1-\frac{r}{\|u_t\|_1}\right)+\frac{2tr\langle\varphi,x\rangle|x_n||\varphi_n|}{\|x\|_1\|u_t\|_1(|x_n|+|x_n+t\varphi_n|)}+\frac{t^2r\langle\varphi,x\rangle\varphi_n^2}{\|x\|_1\|u_t\|_1(|x_n|+|x_n+t\varphi_n|)}}{t|\varphi_n|\left(1+\frac{r}{\|u_t\|_1}\right)+r\|x\|_1t|\varphi_n|\left|\frac{2x_n+t|\varphi_n|}{\|x\|_1\|u_t\|_1(|x_n|+|x_n+t\varphi_n|)}\right|}$$

$$=\limsup_{\substack{t\downarrow 0\\ t<a}} \frac{|\varphi_n|\left(1-\frac{r}{\|u_t\|_1}\right)+\frac{2r\langle\varphi,x\rangle|x_n|}{\|x\|_1\|u_t\|_1(|x_n|+|x_n+t\varphi_n|)}+\frac{tr\langle\varphi,x\rangle|\varphi_n|}{\|x\|_1\|u_t\|_1(|x_n|+|x_n+t\varphi_n|)}}{\left(1+\frac{r}{\|u_t\|_1}\right)+r\left|\frac{2x_n+t|\varphi_n|}{\|u_t\|_1(|x_n|+|x_n+t\varphi_n|)}\right|}$$

$$=\frac{|\varphi_n|\left(1-\frac{r}{\|x\|_1}\right)+\frac{r\langle\varphi,x\rangle}{\|x\|_1^2}}{1+\frac{2r}{\|x\|_1}}$$

$$> 0. \tag{4.8}$$

Case 3. Suppose $\langle\varphi,\ x\rangle < 0$. We can choose *n* such that $\varphi_n \neq 0$ and $x_n\varphi_n < 0$. Then, by $x_n\varphi_n < 0$, we calculate the limit (4.7).

$$\limsup_{\substack{t\downarrow 0\\ t<a}} \frac{t\varphi_n^2\left(1-\frac{r}{\|u_t\|_1}\right)-r\langle\varphi,x\rangle\frac{-2tx_n\varphi_n-t^2\varphi_n^2}{\|x\|_1\|u_t\|_1(|x_n|+|x_n+t\varphi_n|)}}{t|\varphi_n|\left(1+\frac{r}{\|u_t\|_1}\right)+r\|x\|_1\left|\frac{-2tx_n\varphi_n-t^2\varphi_n^2}{\|x\|_1\|u_t\|_1(|x_n|+|x_n+t\varphi_n|)}\right|}$$

$$= \limsup_{\substack{t\downarrow 0\\ t<a}} \frac{t\varphi_n^2\left(1-\frac{r}{\|u_t\|_1}\right)+\frac{2tr\langle\varphi,x\rangle x_n\varphi_n}{\|x\|_1\|u_t\|_1(|x_n|+|x_n+t\varphi_n|)}+\frac{t^2r\langle\varphi,x\rangle\varphi_n^2}{\|x\|_1\|u_t\|_1(|x_n|+|x_n+t\varphi_n|)}}{t|\varphi_n|\left(1+\frac{r}{\|u_t\|_1}\right)+r\|x\|_1 t|\varphi_n|\left|\frac{2x_n+t|\varphi_n|}{\|x\|_1\|u_t\|_1(|x_n|+|x_n+t\varphi_n|)}\right|}$$

$$= \limsup_{\substack{t\downarrow 0\\ t<a}} \frac{t\varphi_n^2\left(1-\frac{r}{\|u_t\|_1}\right)+\frac{2tr|\langle\varphi,x\rangle x_n\varphi_n|}{\|x\|_1\|u_t\|_1(|x_n|+|x_n+t\varphi_n|)}+\frac{t^2r\langle\varphi,x\rangle\varphi_n^2}{\|x\|_1\|u_t\|_1(|x_n|+|x_n+t\varphi_n|)}}{t|\varphi_n|\left(1+\frac{r}{\|u_t\|_1}\right)+r\|x\|_1 t|\varphi_n|\left|\frac{2x_n+t|\varphi_n|}{\|x\|_1\|u_t\|_1(|x_n|+|x_n+t\varphi_n|)}\right|}$$

$$= \limsup_{\substack{t\downarrow 0\\ t<a}} \frac{t\varphi_n^2\left(1-\frac{r}{\|u_t\|_1}\right)+\frac{2tr|\langle\varphi,x\rangle||x_n||\varphi_n|}{\|x\|_1\|u_t\|_1(|x_n|+|x_n+t\varphi_n|)}+\frac{t^2r\langle\varphi,x\rangle\varphi_n^2}{\|x\|_1\|u_t\|_1(|x_n|+|x_n+t\varphi_n|)}}{t|\varphi_n|\left(1+\frac{r}{\|u_t\|_1}\right)+rt|\varphi_n|\left|\frac{2x_n+t|\varphi_n|}{\|u_t\|_1(|x_n|+|x_n+t\varphi_n|)}\right|}$$

$$= \limsup_{\substack{t\downarrow 0\\ t<a}} \frac{|\varphi_n|\left(1-\frac{r}{\|u_t\|_1}\right)+\frac{2r|\langle\varphi,x\rangle||x_n|}{\|x\|_1\|u_t\|_1(|x_n|+|x_n+t\varphi_n|)}+\frac{tr\langle\varphi,x\rangle|\varphi_n|}{\|x\|_1\|u_t\|_1(|x_n|+|x_n+t\varphi_n|)}}{\left(1+\frac{r}{\|u_t\|_1}\right)+r\left|\frac{2x_n+t|\varphi_n|}{\|u_t\|_1(|x_n|+|x_n+t\varphi_n|)}\right|}$$

$$= \frac{|\varphi_n|\left(1-\frac{r}{\|x\|_1}\right)+\frac{r|\langle\varphi,x\rangle|}{\|x\|_1^2}}{1+\frac{2r}{\|x\|_1}}$$

$$> 0. \tag{4.9}$$

(4.6), (4.8) and (4.9) together imply

$$\varphi \notin \widehat{D}^*P_{r\mathbb{B}}\left(x, \frac{r}{\|x\|_1}x\right)(\varphi), \text{ for any } \varphi \in l_\infty \text{ with } \varphi \neq \theta^*.$$

This proves (4.4). Then, part (ii) of this theorem follows from (4.3) and (4.4). □

### 4.5. Banach space $C[0, 1]$

We recall some concepts and notations about $C[0, 1]$, which are used in [15]. Let $(C[0, 1], \|\cdot\|)$ be the Banach space of all continuous real valued functions on [0, 1] with respect to the standard Borel $\sigma$-field $\Sigma$ and with the maximum norm

$$\|f\| = \max_{0\le t\le 1}|f(t)|, \text{ for any } f \in C[0, 1].$$

The dual space of $C[0, 1]$ is denoted by $C^*[0, 1]$ that is $rca[0, 1]$, in which the considered σ-field is the standard Borel σ-field Σ on [0, 1] including all closed and open subsets of [0, 1]. Let $\langle\cdot, \cdot\rangle$ denote the real canonical pairing between $C^*[0,1]$ and $C[0, 1]$. For any $\varphi \in C^*[0,1]$, there is a real valued, regular and countable additive functional $\mu \in rca[0, 1]$, which is defined on the given σ-field Σ on [0, 1], such that

$$\langle\varphi, f\rangle = \int_0^1 f(t)\mu(dt), \text{ for any } f \in C[0, 1].$$

The norm on $C^*[0,1]$ is denoted by $\|\cdot\|_*$ satisfying that, for any $\varphi \in C^*$, the norm of φ in $C^*$ is

$$\|\varphi\|_* = \|\mu\|_* \coloneqq v(\mu, [0, 1]).$$

Where, $v(\mu, [0, 1])$ is the total variation of $\mu$ on [0, 1], which is defined by (see page 160 in [5])

$$v(\mu, [0, 1]) = \sup_{E\in\Sigma}|\mu(E)|.$$

Throughout this subsection, without any special mention, we shall identify φ and $\mu$. The origin of $C[0, 1]$ is denoted by $\theta$, which is the constant function defined on [0, 1] with value 0. The origin of the dual space $C^*$is denoted by $\theta^*$, which is also a constant functional on Σ with value 0. This is,

$$\theta^*(E) = 0, \text{ for every } E \in \Sigma.$$

For any given nonnegative integer $n$, let $\mathcal{P}_n$ denote the closed subspace of $C[0, 1]$ that consists of all real coefficient polynomials of degree less than or equal to $n$. Let $P_{\mathcal{P}_n}$ denote the metric projection from $C[0, 1]$ to $\mathcal{P}_n$. Let $f \in C[0, 1]$ and $p \in \mathcal{P}_n$. Then $p \in P_{\mathcal{P}_n}(f)$ if and only if

$$\begin{aligned}
&\|f - p\| \\
&= \max_{0\le t\le 1}|f(t) - p(t)| \\
&= \min\{\|f - q\|: q(t) = \textstyle\sum_{k=0}^{n} c_k t^k \in \mathcal{P}_n\} \\
&= \min\{\max_{0\le t\le 1}|f(t) - \textstyle\sum_{k=0}^{n} c_k t^k|: \sum_{k=0}^{n} c_k t^k \in \mathcal{P}_n\}.
\end{aligned}$$

The metric projection operator $P_{\mathcal{P}_n}$ from $C[0, 1]$ to $\mathcal{P}_n$ has the following properties proved in [15]. The proofs of these properties are based on the Chebyshev's Equioscillation Theorem (see [4, 15, 25] for details).

**Proposition 5.2 in [15]**. *Let n be a nonnegative integer. Then*

(a) $P_{\mathcal{P}_n}(p) = p$, *for any* $p \in \mathcal{P}_n$;
(b) $P_{\mathcal{P}_n}(f) \neq \emptyset$, *for any* $f \in C[0, 1]$.

**Theorem 5.4 in [15]**. *Let n be a nonnegative integer. The metric projection* $P_{\mathcal{P}_n}: C[0, 1] \to \mathcal{P}_n$ *is a single-valued mapping.*

**Proposition 5.5 in [15]**. *Let n be a nonnegative integer. Let* $f \in C[0, 1]$ *with* $p = P_{\mathcal{P}_n}(f)$. *Then,*

(i) $$\beta p = P_{\mathcal{P}_n}(\beta f), \text{ for any } \beta \in \mathbb{R},$$

$$A(\beta f, \beta p) = |\beta| A(f, p) \quad \text{and} \quad S(\beta f, \beta p) = S(f, p),$$

$$\epsilon(\beta f, \beta p) = \epsilon(f, p), \text{ if } \beta \ge 0 \quad \text{and} \quad \epsilon(\beta f, \beta p) = -\epsilon(f, p), \text{ if } \beta < 0.$$

(ii) $$p + q = P_{\mathcal{P}_n}(f + q), \text{ for any } q \in \mathcal{P}_n,$$

$$A(f + q, p + q) = A(f, p) \quad \text{and} \quad S(f + q, p + q) = S(f, p).$$

*In particular, we have*

$$p + \lambda = P_{\mathcal{P}_n}(f + \lambda), \quad \textit{for any } \lambda \in \mathbb{R}.$$

**Lemma 4.5**. *For any positive integer n, define* $A_n \in C[0, 1]$ *by the following* $(n + 1) \times (n + 1)$ *determinate*

$$A_n(t) = \begin{vmatrix} 1 & 1 & 1 & 1 & \cdots & 1 & 1 \\ 1 & t & t^2 & t^3 & \cdots & t^{n-1} & t^n \\ 1 & t^2 & t^4 & t^6 & \cdots & t^{2(n-1)} & t^{2n} \\ 1 & t^3 & t^6 & t^9 & \vdots & t^{3(n-1)} & t^{3n} \\ \cdots & \cdots & \cdots & \cdots & & \cdots & \cdots \\ 1 & t^{n-1} & t^{2(n-1)} & t^{3(n-1)} & \cdots & t^{(n-1)^2} & t^{(n-1)n} \\ 1 & t^n & t^{2n} & t^{3n} & \cdots & t^{(n-1)n} & t^{n^2} \end{vmatrix}, \quad \textit{for } t \in [0, 1].$$

*Then,*

$$A_n(t) \neq 0, \textit{ for every } t \in (0, 1). \tag{4.10}$$

*Proof*. We prove (4.10) by the method of induction.

Step 1. For $n = 2$. We have

$$A_2(t) = \begin{vmatrix} 1 & 1 & 1 \\ 1 & t & t^2 \\ 1 & t^2 & t^4 \end{vmatrix} = t^5 - 2t^4 + 2t^2 - t < 0, \text{ for every } t \in (0, 1).$$

Step 2. Assume that

$$A_{n-1}(t) \neq 0, \text{ for every } t \in (0, 1).$$

Step 3. We prove (4.10) for $n$. For any $t \in (0, 1)$, By the assumption in step 2, we have

$$A_n(t) = \begin{vmatrix} 1 & 1 & 1 & 1 & \cdots & 1 & 1 \\ 1 & t & t^2 & t^3 & \cdots & t^{n-1} & t^n \\ 1 & t^2 & t^4 & t^6 & \cdots & t^{2(n-1)} & t^{2n} \\ 1 & t^3 & t^6 & t^9 & \vdots & t^{3(n-1)} & t^{3n} \\ \cdots & \cdots & \cdots & \cdots & & \cdots & \cdots \\ 1 & t^{n-1} & t^{2(n-1)} & t^{3(n-1)} & \cdots & t^{(n-1)^2} & t^{(n-1)n} \\ 1 & t^n & t^{2n} & t^{3n} & \cdots & t^{(n-1)n} & t^{n^2} \end{vmatrix}$$

$$= \begin{vmatrix} 1 & 1 & 1 & 1 & \cdots & 1 & 1 \\ 0 & t-1 & t^2-1 & t^3-1 & \cdots & t^{n-1}-1 & t^n-1 \\ 0 & t^2-1 & t^4-1 & t^6-1 & \cdots & t^{2(n-1)}-1 & t^{2n}-1 \\ 0 & t^3-1 & t^6-1 & t^9-1 & \vdots & t^{3(n-1)}-1 & t^{3n}-1 \\ \cdots & \cdots & \cdots & \cdots & & \cdots & \cdots \\ 0 & t^{n-1}-1 & t^{2(n-1)}-1 & t^{3(n-1)}-1 & \cdots & t^{(n-1)^2}-1 & t^{(n-1)n}-1 \\ 0 & t^n-1 & t^{2n}-1 & t^{3n}-1 & \cdots & t^{n(n-1)}-1 & t^{n^2}-1 \end{vmatrix}$$

$$
= \begin{vmatrix}
1 & 0 & 0 & 0 & \cdots & 0 & 0 \\
0 & t-1 & t^2-1 & t^3-1 & \cdots & t^{n-1}-1 & t^n-1 \\
0 & t^2-1 & t^4-1 & t^6-1 & \cdots & t^{2(n-1)}-1 & t^{2n}-1 \\
0 & t^3-1 & t^6-1 & t^9-1 & \vdots & t^{3(n-1)}-1 & t^{3n}-1 \\
\cdots & \cdots & \cdots & \cdots & & \cdots & \cdots \\
0 & t^{n-1}-1 & t^{2(n-1)}-1 & t^{3(n-1)}-1 & \cdots & t^{(n-1)^2}-1 & t^{(n-1)n}-1 \\
0 & t^n-1 & t^{2n}-1 & t^{3n}-1 & \cdots & t^{n(n-1)}-1 & t^{n^2}-1
\end{vmatrix} =
$$

$$
\begin{vmatrix}
1 & 0 & 0 & 0 & \cdots & 0 & 0 \\
0 & t-1 & t^2-1 & t^3-1 & \cdots & t^{n-1}-1 & t^n-1 \\
0 & t^2-t & t^4-t^2 & t^6-t^3 & \cdots & t^{2(n-1)}-t^{n-1} & t^{2n}-t^n \\
0 & t^3-t^2 & t^6-t^4 & t^9-t^6 & \vdots & t^{3(n-1)}-t^{2(n-1)} & t^{3n}-t^{2n} \\
\cdots & \cdots & \cdots & \cdots & & \cdots & \cdots \\
0 & t^{n-1}-t^{n-2} & t^{2(n-1)}-t^{2(n-2)} & t^{3(n-1)}-t^{3(n-2)} & \cdots & t^{(n-1)^2}-t^{(n-1)(n-2)} & t^{(n-1)n}-t^{(n-2)n} \\
0 & t^n-t^{n-1} & t^{2n}-t^{2(n-1)} & t^{3n}-t^{3(n-1)} & \cdots & t^{n(n-1)}-t^{(n-1)^2} & t^{n^2}-t^{(n-1)n}
\end{vmatrix}
$$

$$
= (t-1)(t^2-1)(t^3-1)\ldots(t^{n-1}-1)(t^n-1)
$$

$$
\begin{vmatrix}
1 & 0 & 0 & 0 & \cdots & 0 & 0 \\
0 & 1 & 1 & 1 & \cdots & 1 & 1 \\
0 & t & t^2 & t^3 & \cdots & t^{n-1} & t^n \\
0 & t^2 & t^4 & t^6 & \vdots & t^{2(n-1)} & t^{2n} \\
\cdots & \cdots & \cdots & \cdots & & \cdots & \cdots \\
0 & t^{n-2} & t^{2(n-2)} & t^{3(n-2)} & \cdots & t^{(n-2)(n-1)} & t^{(n-2)n} \\
0 & t^{n-1} & t^{2(n-1)} & t^{3(n-1)} & \cdots & t^{(n-1)^2} & t^{(n-1)n}
\end{vmatrix}
$$

$$
= (t-1)(t^2-1)(t^3-1)\ldots(t^{n-1}-1)(t^n-1)\, t\, t^2 \ldots t^{n-2} t^{n-1}
$$

$$
\begin{vmatrix}
1 & 0 & 0 & 0 & \cdots & 0 & 0 \\
0 & 1 & 1 & 1 & \cdots & 1 & 1 \\
0 & 1 & t & t^2 & \cdots & t^{n-2} & t^{n-1} \\
0 & 1 & t^2 & t^4 & \vdots & t^{2(n-2)} & t^{2(n-1)} \\
\cdots & \cdots & \cdots & \cdots & & \cdots & \cdots \\
0 & 1 & t^{n-2} & t^{2(n-2)} & \cdots & t^{(n-2)^2} & t^{(n-2)(n-1)} \\
0 & 1 & t^{n-1} & t^{2(n-1)} & \cdots & t^{(n-2)(n-1)} & t^{(n-1)^2}
\end{vmatrix}
$$

$$
= (t-1)(t^2-1)(t^3-1)\ldots(t^{n-1}-1)(t^n-1)\, t\, t^2 \ldots t^{n-2} t^{n-1} A_{n-1}(t)
$$

$\neq 0$, for every $t \in (0, 1)$.

This proves (4.10). □

Let $D$ be a nonempty subset of $\mathcal{P}_n$. Let Coe($D$) be the collection of all coefficients of polynomials in $D$, which is called the coefficient set of $D$. That is,

$$
\mathrm{Coe}(D) = \{a \in \mathbb{R}\text{: there is } p \in D \text{ such that } a \text{ is a coefficient of } p\}.
$$

**Proposition 4.6**. *Let n be a positive integer. Let D be a nonempty subset of* $\mathcal{P}_n$ *with coefficient set* Coe(*D*). *Then, we have*

$$D \text{ is } \|\cdot\|\text{-bounded} \quad \Leftrightarrow \quad \text{Coe}(D) \text{ is } |\cdot|\text{-bounded}.$$

*Proof.* "⇐". Suppose that Coe(*D*) is |·|-bounded by a positive number *d*. Since we only consider the polynomials on [0, 1], this implies that *D* is ||·||-bounded by $(n+1)d$.

"⟹". Suppose that *D* is ||·||-bounded by *b*, for some positive number *b*. Take an arbitrary element $p \in D$ with

$$p(t) = a_0 + a_1 t + a_2 t^2 + \dots + + a_n t^n, \text{ for } t \in [0, 1]. \tag{4.11}$$

Respectively substituting *t* by 1, $\frac{1}{2}$, $\frac{1}{2^2}$, …, $\frac{1}{2^n}$ in (4.11), we calculate the corresponding values and obtain a system of linear equations.

$$\begin{aligned} &a_0 + a_1 + a_2 + \dots + a_n = b_0, \\ &a_0 + a_1\frac{1}{2} + a_2\frac{1}{2^2} + \dots + a_n\frac{1}{2^n} = b_1, \\ &a_0 + a_1\frac{1}{2^2} + a_2\frac{1}{2^4} + \dots + a_n\frac{1}{2^{2n}} = b_2, \\ &\;\dots\dots \\ &a_0 + a_1\frac{1}{2^n} + a_2\frac{1}{2^{2n}} + \dots + a_n\frac{1}{2^{n^2}} = b_n. \end{aligned} \tag{4.12}$$

Where $b_0, b_1, b_2, \dots, b_n$ are real numbers. By the assumption that *D* is ||·||-bounded by *b*, we have

$$|b_k| \le b, \text{ for } k = 0, 1, 2, \dots, n.$$

By the definition of $A_n(\cdot)$, the system of linear equations (4.12) in terms of the "variables" $a_0$, $a_1$, $a_2$, …, $a_n$ can be rewritten as

$$A_n\left(\frac{1}{2}\right)\begin{pmatrix} a_0 \\ a_1 \\ a_2 \\ \cdots \\ a_n \end{pmatrix} = \begin{pmatrix} b_0 \\ b_1 \\ b_2 \\ \cdots \\ b_n \end{pmatrix}. \tag{4.12}$$

For $k = 0, 1, 2, \dots, n$, let $A_n^k\left(\frac{1}{2}\right)$ denote the $(n+1)\times(n+1)$ determinate that is defined by substituting the *k*th-column in $A_n\left(\frac{1}{2}\right)$ by $(b_0, b_1, b_2, \dots, b_n)^T$. By Lemma 4.5, $A_n\left(\frac{1}{2}\right) \neq 0$, so (4.12) has a unique solution for $a_0, a_1, a_2, \dots, a_n$ in term of $b_0, b_1, b_2, \dots, b_n$. We have

$$|a_k| = \frac{\left|A_n^k\left(\frac{1}{2}\right)\right|}{\left|A_n\left(\frac{1}{2}\right)\right|} \le \frac{n!b}{\left|A_n\left(\frac{1}{2}\right)\right|}, \text{ for } k = 0, 1, 2, \dots, n. \tag{4.13}$$

This implies that Coe(*D*) is |·|-bounded by $\frac{n!b}{\left|A_n\left(\frac{1}{2}\right)\right|}$, which does not depend on $p \in D$. □

**Corollary 4.7**. *Let n be a positive integer. Let D be a nonempty subset of* $\mathcal{P}_n$ *with coefficient set* Coe(*D*). *Define* $D' = \{p' \in \mathcal{P}_{n-1} : p \in D\}$. *Then, we have*

$$D \text{ is } \|\cdot\|\text{-bounded} \quad \Longrightarrow \quad D' \text{ is } \|\cdot\|\text{-bounded}.$$

*Proof*. Suppose that $D$ is $\|\cdot\|$-bounded by $b$, for some $b > 0$. By (4.13) in the proof of Proposition 4.6, Coe($D$) is $|\cdot|$-bounded by $\frac{n!b}{\left|A_n\left(\frac{1}{2}\right)\right|}$. This implies that Coe($D'$) is $|\cdot|$-bounded by $\frac{n!nb}{\left|A_n\left(\frac{1}{2}\right)\right|}$. By Proposition 4.6 again, $D'$ is $\|\cdot\|$-bounded by $\frac{n!n^2b}{\left|A_n\left(\frac{1}{2}\right)\right|}$. □

Let $n$ be a positive integer. In [15], it is proved that $\mathcal{P}_n$ is closed and not open in $C[0, 1]$ and the metric projection $P_{\mathcal{P}_n}: C[0, 1] \to \mathcal{P}_n$ is a single valued mapping. In the following theorem, we prove the continuity of $P_{\mathcal{P}_n}$.

**Theorem 4.8**. *Let $n$ be a positive integer. The metric projection $P_{\mathcal{P}_n}: C[0, 1] \to \mathcal{P}_n$ is a continuous mapping.*

*Proof*. Let $\{f_m\}$ be a convergent sequence in $C[0, 1]$ with limit $f$. Then $\{f_m\}$ is $\|\cdot\|$-bounded. Let $P_{\mathcal{P}_n}(f_m) = p_m$, for all $m$. Take an arbitrary $w \in \mathcal{P}_n$, we have

$$\|p_m\| \leq \|p_m - f_m\| + \|f_m\| \leq \|w - f_m\| + \|f_m\| \leq \|w\| + 2\|f_m\|.$$

This implies that $\{p_m\}$ is $\|\cdot\|$-bounded. By Corollary 4.7, $\{p_m'\}$ is $\|\cdot\|$-bounded. We define the upper bound $d = \max\{\|p_m'\|: m = 1, 2, 3 \dots\}$. For any given $\varepsilon > 0$, we take $\delta > 0$ with $\delta < \frac{\varepsilon}{d}$. By the mean value theorem, this implies that, for any $s, t \in [0, 1]$, we have

$$|s - t| < \delta \;\; \Longrightarrow \;\; |p_m(s) - p_m(t)| < \varepsilon, \text{ for every } m = 1, 2, \dots .$$

This implies that $\{p_m\}$ is a bounded and equicontinuous subset in $C[0, 1]$. By Azsela-Ascali Theorem, $\{p_m\}$ is conditionally compact (see [5]). Since $\mathcal{P}_n$ is closed in $C[0, 1]$, then, there is a subsequence $\{p_{m_i}\} \subseteq \{p_m\} \subseteq \mathcal{P}_n$ and there is $p \in \mathcal{P}_n$ such that $p_{m_i} \to p$, as $i \to \infty$. For any $q \in \mathcal{P}_n$, by $P_{\mathcal{P}_n}(f_m) = p_m$, we have

$$\begin{aligned} \|f - p\| &\leq \|f - f_{m_i}\| + \|f_{m_i} - p_{m_i}\| + \|p - p_{m_i}\| \\ &\leq \|f - f_{m_i}\| + \|f_{m_i} - q\| + \|p - p_{m_i}\| \\ &\leq \|f - f_{m_i}\| + \|f_{m_i} - f\| + \|f - q\| + \|p - p_{m_i}\|. \end{aligned}$$

By $p_{m_i} \to p$ and $f_{m_i} \to f$ as $i \to \infty$, this implies

$$\|f - p\| \leq \|f - q\|, \text{ for any } q \in \mathcal{P}_n.$$

It follows that $P_{\mathcal{P}_n}(f) = p$. This proves the continuity of $P_{\mathcal{P}_n}$. □

Let $n$ be a positive integer. By Theorem 5.4 in [15] and Theorem 4.8, the metric projection $P_{\mathcal{P}_n}: C[0, 1] \to \mathcal{P}_n$ is a single valued continuous mapping. Then, by Lemma 3.1, we have the following results for the calculation of fixed points of $\widehat{D}^*P_{\mathcal{P}_n}(\cdot)$.

**Corollary 4.10**. *Let $n$ be a positive integer. The fixed-points of Mordukhovich differential operator $\widehat{D}^*P_{\mathcal{P}_n}(\cdot)$ of the metric projection $P_{\mathcal{P}_n}: C[0, 1] \to \mathcal{P}_n$ can be calculated as follows. For*

$h \in C[0,1]$, *we have*

$$\mathcal{F}\left(\widehat{D}^* P_{\mathcal{P}_n}(h)\right) = \left\{y^* \in C^*[0,1]: \limsup_{g \to h} \frac{\langle y^*,\ (g-h)-(P_{\mathcal{P}_n}(g)-P_{\mathcal{P}_n}(h))\rangle}{\|g-h\|+\|P_{\mathcal{P}_n}(g)-P_{\mathcal{P}_n}(h)\|} \leq 0\right\}.$$

We write the orthogonal space of $\mathcal{P}_n$ by $\mathcal{P}_n^{\perp}$, which is defined as follows.

$$\mathcal{P}_n^{\perp} = \{\mu \in\ C^*[0,1]: \langle \mu,\ p\rangle = 0, \text{for all } p \in \mathcal{P}_n\}\ .$$

**Theorem 4.11**. *Let n be a positive integer. Let* $f \in C[0,\ 1]$, $p = P_{\mathcal{P}_n}(f)$ *and* $\mu \in C^*[0,\ 1]$ *with* $\langle \mu, f\rangle \neq\ 0$, *we have*

(i) $\mu \notin \widehat{D}^* P_{\mathcal{P}_n}(f,p)(\gamma)$, *for any* $\gamma \in \mathcal{P}_n^{\perp}$;

(ii) $\mu \in \mathcal{P}_n^{\perp} \Longrightarrow\ \mu \notin \mathcal{F}\left(\widehat{D}^* P_{\mathcal{P}_n}(f,p)\right)$.

*Proof.* Proof of (i). For the given $f \in C[0,\ 1]$ and $p = P_{\mathcal{P}_n}(f)$, let $\mu \in C^*[0,1]$ with $\langle \mu, f\rangle \neq 0$. Without loose of the generality, we assume that $\langle \mu, f\rangle > 0$. Let $\gamma \in \mathcal{P}_n^{\perp}$. We take a line segment direction as $g_\lambda = f + \lambda f$, for $\lambda \downarrow 0$. One has

$$g_\lambda = f + \lambda f\ \to f, \text{ as } \lambda \downarrow 0.$$

By Proposition 5.5 in [15], we have that $P_{\mathcal{P}_n}(g_\lambda) = P_{\mathcal{P}_n}(f + \lambda f) = (1+ \lambda)\, p$. This implies that $P_{\mathcal{P}_n}(g_\lambda) \to p$, as $\lambda \downarrow 0$. By $\gamma \in \mathcal{P}_n^{\perp}$, we calculate

$$\limsup_{\substack{(g,q)\to(f,p)\\ q=P_{\mathcal{P}_n}(g)}} \frac{\langle \mu, g-f\rangle - \langle \gamma,\ q-p\rangle}{\|g-f\| + \|q-p\|}$$

$$\geq \limsup_{(g_\lambda, P_{\mathcal{P}_n}(g_\lambda))\to(f,p)} \frac{\langle \mu, g_\lambda - f\rangle - \langle \gamma,\ P_{\mathcal{P}_n}(g_\lambda) - p\rangle}{\|g_\lambda - f\| + \|P_{\mathcal{P}_n}(g_\lambda) - p\|}$$

$$= \limsup_{(f+\lambda h, P_{\mathcal{P}_n}(f+\lambda h))\to(f,p)} \frac{\lambda\langle \mu, f\rangle}{\lambda\|f\| + \lambda\|p\|}$$

$$= \frac{\langle \mu, f\rangle}{\|f\| + \|p\|} > 0.$$

This implies that $\mu \notin \widehat{D}^* P_{\mathcal{P}_n}(f,p)(\gamma)$. Then, part (ii) follows from part (i) immediately. □

**Acknowledgments.** The author is very grateful to Professor Akhtar A. Khan, Professor Boris S. Mordukhovich, Professor Simeon Reich, Professor Christiane Tammer and Professor Jen-Chih Yao for their kind communications, valuable suggestions and enthusiasm encouragements in the development stage of this work.